\documentclass[12pt]{article}
\usepackage{mathptmx}
\usepackage{graphicx}
\usepackage{float}
\usepackage{bm}
\usepackage{amssymb}
\makeatletter
\author{Robert P. C. de Marrais \footnote{Email address:  rdemarrais@alum.mit.edu} \\ \noindent
\emph{Thothic Technology Partners, P.O.Box 3083, Plymouth MA 02361}}

\title{Voyage by Catamaran:\linebreak  Effecting Semantic Network \textit{Bricolage} via Infinite-Dimensional Zero-Divisor Ensembles}

\begin{document}
\maketitle
\makeatother

\begin{abstract}
Continuing arguments presented [1] or announced [2][3] in
\textit{Complex Systems}, zero-divisor (ZD) foundations for
scale-free networks (evinced, in particular, in the ``fractality''
of the Internet) are decentralized.  \textit{Spandrels} -- quartets
of ZD-free or ``hidden'' box-kite-like structures (HBK's) in the
$2^{N+1}$-ions -- are ``exploded'' from (and uniquely linked to)
each standard box-kite in the $2^{N}$-ions, $N \geq 4$.  Any HBK
houses, in a \textit{cowbird's nest}, exactly one copy of the
(ZD-free) octonions -- the recursive basis for all ZD ensembles, and
each a potential waystation for alien-ensemble infiltration in the
large, or metaphor-like jumps in the small.  Cowbirding models what
evolutionary biologists[4], and structural mythologist Claude
L\'evi-Strauss before them[5], term \textit{bricolage}: the
opportunistic co-opting of objects designed for one purpose to serve
others unrelated to it. Such arguments entail a switch of focus,
from the octahedral box-kite's four triangular \textit{sails}, to
its trio of square \textit{catamarans} and their box-kite-switching
\textit{twist products}.
\end{abstract}

\section{From Box-Kites to Brocades via Catamaran Twists}
This work had its beginnings in [6], where an abstract result of
Guillermo Moreno [7] was employed to explicitly delineate the ZD
structure of the 16-D Sedenions. These hypercomplex numbers are
reached via the Cayley-Dickson Process (CDP), a dimension-doubling
algorithm which begins with the linear Real Numbers, moves to the
Complex Plane, generates the quaternions' non-commutative 4-space,
then the 8-D non-associative octonions:  all so many waystations
\textit{en route} to the ``pathology'' of zero-divisors, found in
all $2^{N}$-ions, $N \geq 4$.

The key to the early results was found in simplifying CDP itself,
reducing it to a set of two bit-twiddling rules, exploiting one
convention. The quaternions' imaginary units can be represented in
two different ways with the subscripts 1, 2, 3. (Which way obtains
depends upon which order of multiplication of two of them yields a
positively signed instance of the third). The octonions' 7 non-real
units, though, can be indexed in 480 distinct ways, perhaps a dozen
of which are in actual use among various specialists.  The
sedenions' indexing schemes number in the billions.  Yet only one
kind of scheme can work for \textit{all} $2^{N}$-ions: first, index
their units with the integers $1$ through $2^{N} - 1$, with $0$
denoting the Real unit; then, assume the index of a unit produced by
multiplying two others is the XOR of their indices.

Note down only the indices, suppressing the usual tedium of writing
a lower-case $i$ with a hard-to-read subscript, and list XOR sets in
parenthesized triplets, with the first two units ordered so that
their product, in the third slot, has positive sign.  By such
\textit{cyclical positive ordering} (CPO), the two possible
quaternion multiplication rules are written $(1, 2, 3)$ and $(1, 3,
2)$, with the former taken here as basis for all CDP recursion. And,
the 480 octonion rule-sets collapse to one set of 7 CPO triplets --
\textit{trips} for short -- corresponding precisely and only to the
7 associative triplets in this otherwise nonassociative number
space: $(1, 2, 3)$; $(1, 4, 5)$; $(1, 7, 6)$; $(2, 4, 6)$; $(2, 5,
7)$; $(3, 4, 7)$; $(3, 6, 5)$.  All these remain associative
triplets in all higher $2^{N}$-ions -- a fact we call \textit{Rule
0}. And, with the index-0 Real unit appended to the set, each also
provides a true copy of the quaternions.

The next crucial step is to understand how and why some trips thus
derived are \textit{not} in \textit{counting order}.  Assume, as
standard CDP does, that the $2^{N} - 1$ units of a given set of
$2^{N}$-ions can be multiplied on the right by a new unit whose
index exceeds all theirs, the \textit{generator} {\bf G} of the
$2^{N + 1}$-ions -- to yield resultant units with new and higher
indices, all with positive sign.  {\bf G} is just the unique unit of
index $2^{N}$, and for each unit with index $L < {\bf G}$, the CPO
trip can be written like this:  $( L, {\bf G}, {\bf G} + L)$. (Note
that the product's index is identical to the XOR ${\bf G} \veebar
L$, since {\bf G}, by definition, can be represented by a bit to the
left of any L's bitstring expression.) This is \textit{Rule 1}. It
completely explains the indexing of the quaternions: if the usual
imaginary unit has index $1$, then {\bf G}$ = 2$ and Rule 1 yields
$(1, 2, 3)$.  For the octonions, we inherit the quaternion's
singleton trip, and generate 3 more with Rule 1:  $(1, 4, 5); (2, 4,
6); (3, 4, 7)$.  But how do we get the 3 that remain? Here's where
we invoke \textit{Rule 2}.

The quaternions need no Rule 2, so start an induction by assuming it
works only on Rule 0 trips from prior CDP generations. For any such
trip, hold one index fixed, then add {\bf G} to each of the other
two and \textit{switch their positions}. For the octonions, with
only one Rule 0 trip to manipulate, we get the needed three Rule 2
trips by this tactic: fixing $1, 2,$ and $3$ in that order, we get
$(1, 3+4, 2+4); (3+4, 2, 1+4); (2+4, 1+4, 3)$. Cyclical rotation
bringing the smallest index to the left yields the 3 extra trips
written above: $(1, 7, 6); (2, 5, 7); (3, 6, 5)$.  All applications
of CDP to the standard real and imaginary units, for $N$ as large as
desired, are completely covered by these rules. Also, the total
number of trips, which simple combinatorics tell us we can generate
in a given set of $2^{N}$-ions, is just $(2^{N}-1)(2^{N}-2)/3!$ --
hence, $1$ for the quaternions (where $N = 2$); $7$ for the $N = 3$
octonions; $35$ for the sedenions where ZDs are first in evidence;
and $155$ for the 32-D pathions, where the signature of scale-free
behavior, as evidenced in the World Wide Web's implicit
\textit{fractality} (Sir Tim Berners-Lee's term for it [8]), is
first revealed.

A fine point:  fractality has two related but distinct senses. That
of ZD structures concerns the one-to-one mapping one can make
between points in a classic 2-D fractal and the empty cells in our
Emanation Tables (ETs). As derived in [2] and illustrated in [9],
these are spreadsheet-like multiplication tables of zero-divisors
whose unfilled cells indicates row and column ZD entries which do
not mutually zero-divide \textit{each other}. The Web's fractality,
though, concerns statistical distributions -- of links, say, between
web pages or routers, which have much higher densities at some nodes
than others, with density distributions being far less Gaussian and
normal than they tend toward being Mandelbrot-set self-similar. Our
ET cell entries (pairs of which sharing symmetric row and column
labels map directly to pairs of ZD-saturated diagonals, in some
plane associated with some one of some box-kite's six vertices) do
not remain pure Number Theory entities, but become heuristically
placed statistical markers, when dynamic models of actual networks
are simulated and/or searched by means of ZD ensembles.  As
construction of the methodology for such model-building is our
ultimate aim (only partially realized in these pages), we feel
justified in assuming the applicability of fractality, in both its
senses, to the agenda being sketched here. Further elaboration on
this point must be deferred to future studies.

Now, to understand what happens in 32-D, we must first explain ZDs'
workings in the sedenions.  Moreno's abstract treatment of their
interrelationship was framed in the physicist's favored language of
semi-simple Lie groups:  the largest exceptional group, $E_{8}$, has
240 roots which form a loop (the non-associative equivalent of a
group) isomorphic to the unit octonions; the automorphism group of
$E_{8}$, the smallest exceptional group, $G_{2}$, is homomorphic to
the symmetry patterns displayed by ZDs in the sedenions.  And, since
this same $G_{2}$ is also the basis of the derivation algebra that
recursively creates (via CDP) the $2^{N + 1}$-ions from the
$2^{N}$-ions, for all $N > 4$, he would argue this same homomorphism
obtains for all such $N.$  But homo- (as opposed to iso- ) morphism
is a rather imprecise tool for obtaining anything like concrete
results.  The approach taken in [6]:  use minimal assumptions and
bit-twiddling rules.

Since $i^{n} \neq 0$ for any imaginary unit $i$ of any index, raised
to any finite power $n$, the simplest ZD must entail the sum or
difference of a \textit{pair} of imaginaries, and zero will only
result from the product of at least two such pairings. Rather simple
by-hand calculations quickly showed one such unit must have index $L
< {\bf G}$, and its partner have index $U > {\bf G}$ \textit{not the
XOR of $L$ with {\bf G}}. This meant one could pick any octonion (7
choices) and match it with any of the 6 suitable sedenions with
index $> 8$, making for 42 planes or \textit{assessors} whose
diagonal line-pairs contain only (and all the) ZDs.  But these 84
lines do not all mutually zero-divide with each other; those which
\textit{do}, have their behavior summarized in 7 geometrically
identical diagrams, the octahedral wireframe figures called
box-kites.  Their manner of assembly was determined by 3 simple
production rules.

Label the 3 vertices of some triangle among the octahedral grid's 8
with the letters $A$, $B$, $C$, and those of the opposite face $F$,
$E$, $D$, so that these are at opposite ends of lines through the
center {\bf S}  --  $AF$, $BE$, $CD$  --  which we call
\textit{struts}. Assume each vertex represents a plane whose two
units are indicated by the same letter, in upper or lower case
depending on whether the index is greater or less than {\bf G} --
$U$ and $L$ indices respectively. Call {\bf S}, the seventh octonion
index not found on a vertex, the \textit{strut constant}, and use it
to distinguish the 7 box-kites, each of which contains but 6 of the
42 sedenion assessors. For any chosen {\bf S}, there will be 3 pairs
of octonions forming trips with it, and the indices forming such
pairs are placed on \textit{strut-opposite} vertices (i.e., at ends
of the same strut, not edge).  Neither diagonal at one end of a
strut will mutually zero-divide with either at the other:  some $k
\cdot (A \pm a)$ will not yield zero when multiplied by any $q \cdot
(F \pm f)$, $k$ and $q$ arbitrary real scalars.  But either
diagonal, at any assessor, produces zero when multiplied by
\textit{exactly one} of the assessor diagonals at the other end of a
shared edge. Half the edges have ``[+]'' edge-currents (the
diagonals slope the same way, as with $(A + a) \cdot (D + d) = (A -
a) \cdot (D - d) = 0$), while the other six have edges marked
\linebreak ``[-]'' (e.g., $(A + a) \cdot (B - b) = (A - a) \cdot (B
+ b) = 0)$. With these conventions, we can assert the production
rules.

First, if diagonals at A and B mutually zero-divide, each also does
so with a diagonal of C (oppositely signed copies of whose unit
pairings embody the zero produced when A and B diagonal
unit-pairings are multiplied ):  A and B \textit{emanate} C, from
whence the \textit{emanation tables} (ETs) we'll see presently
(where A and B display as row and column labels, designating a
spreadsheet cell with content C). Corollarily, their L-indices
$(\textit{a, b, c})$ form a trip only if their assessors' diagonals
each mutually zero-divide one of those at each of the other two. A
\textit{sail} is such a triad of assessors, representable by a
triangle on the box-kite.  As shown in [1], there is \textit{exactly
one} sail per box-kite with all three edges marked ``[-]'':  the
\textit{zigzag}, so called because its 6-cycle of zero-divisions,
determined by tracing its edges twice, shows an alternation of
\verb|/| and \verb|\| slopings among the diagonals sequentially
engaged in product-forming.  By convention, its assessors are A, B,
and C, with $(\textit{a, b, c})$ in CPO order, rotated to make
$\textit{a}$ the smallest integer.

Third, any assessor belongs to 2 sails, implying 4 in all, touching
only at vertices, like same-color checkerboard squares.  The pair of
imaginary units forming any assessor always split so that one has
index less than the generator G (the \textit{L-}unit, L for lower),
and one has index greater than G (the \textit{U-}unit, U for upper).
An assessor's \textit{L}-unit is written with the same letter as the
assessor proper, but in lowercase italics, while the \textit{U}-unit
is written in uppercase italics:  the pair of units designated
"assessor A" is thus equivalent to \textit{(a, A)}.  The
\textit{L-}units of each sail form associative triplets, hence
L-trips. The \textit{trefoil} L-trips are $(\textit{a, d, e});
(\textit{d, b, f}); (\textit{e, f, c})$ -- with leftmost terms not
necessarily the smallest in their trios, and each derived from the
zigzag L-trip (\textit{Z-trip}) by flipping L-indices along 2
struts, holding \textit{$a$, $b$, $c$} each fixed in turn. The
remaining 4 triangles are \textit{vents}, with the face opposite the
zigzag, DEF, understood as meant when written with a capital ``V.''
The alternation of sails (made of colored paper maybe) with empty
spaces ``where the wind blows,'' and the kite-like structural
stability implied by the 3 ZD-free orthogonal struts (made of wooden
or plastic doweling, perhaps), motivates the conceit of calling
these \textit{box-kites} in the first place. As vent and zigzag use
up all 6 negative edge-currents, edge-currents joining trefoil-based
assessors D, E, F to the zigzag's A, B, C are \textit{positive}.

Previous work focused on sails, whose algebraic closure, and
capacity for recursive construction for growing $N$, make them
exceedingly rich sources of structural information.  But the second
production rule is where our interest will focus here: L- (or U- )
indices can be swapped (with a sign flip) between assessors sharing
an edge, yielding assessor pairs in other box-kites with different
{\bf S}. Hence, since $(A + a) \cdot (B - b) = 0$, then so will $(A
+ b) \cdot (B + a)$ -- with caveats for $N > 4$ if the box-kite is
\textit{Type II}, which we'll soon get to.  Opposite edges of the
same square (one of 3 mutually orthogonal ones) twist to the same
box-kite and have the same edge-sign.  Quandrangular
\textit{catamarans} -- like triangular sails -- have a richness all
their own.  (See Figure 1.)

In the sedenions, all box-kites are \textit{Type I}:  for any zigzag
assessor Z and its Vent strut-opposite V, $(z, {\bf S}, v)$ and $(Z,
{\bf S}, V)$ are CPO. (For Type II's, two of the strut's trips have
reversed orientations.) For all Type I's of any $N$, all 3
catamarans share an invariant feature: the orientation of L-trip
products along each edge is counterclockwise along 3 successive
sides, with the fourth (with negative edge-sign) showing a clockwise
reversal.

\begin{figure}[H]
\centering
\includegraphics[scale = .4] {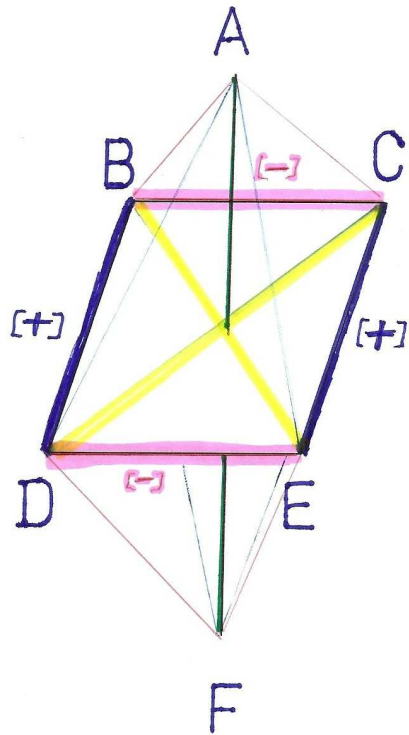}
\caption{Parallel edges of catamarans (one perpendicular to each
strut in an octahedral box-kite) \textit{twist} into assessor pairs
with oppositely signed edge-currents, in a box-kite with different
\textit{strut constant}:  BC and DE, both in sails completed by A,
have twist products with {\bf S} = $f$; for DB and CE, completed by
F, twistings have {\bf S} = $a$.  The $5^{th}$ and $6^{th}$
(necessarily strut-opposite) assessors in each are found by twisting
$(A, a)$ and $(F, f)$ with $({\bf X, S})$ -- assumed at the center,
where struts intersect -- double-covering mast and keel
respectively.}
\end{figure}

Catamarans orthogonal to struts AF, BE, CD have reversed edges DE,
FD, EF, respectively. Rotate their frames to put the reversed edge
on top, and shade or color such edges to specify their catamaran.
Then draw two more catamaran boxes, to the right and just below
this. The top and bottom edges on the right display L- and U- index
twists from the starting box respectively; the left and right edges
below show L and U twists from their vertical counterparts above.
Figure 2 is an instance of such a \textit{Royal Hunt Diagram}, after
the $5^{th}$ moving-line text of \textit{I Ching} Hexagram 8,
Holding Together: ``In the hunt, the king uses beaters on three
sides only and foregoes game that runs off the front.''

\begin{figure}[H]
\centering
\includegraphics[width = 4in, height = !] {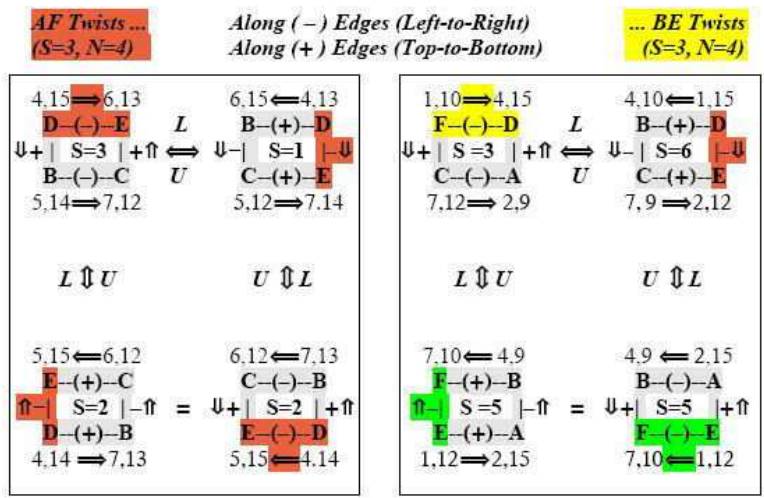}
\caption{Royal Hunt Diagram.  The 3 catamaran {\bf S}'s form an
associative triple:  the bottom-left box, twisted a second time
along its other set of parallels, yields the same resultant as the
second twist of the top-right. (The two bottom boxes differ only by
a $90^{o}$ rotation.) Twists involving HBK's and Type II box-kites
also show a double twist of another kind: the proper tracing order
along the perimeter of the box twisted \textit{to} will be
\textit{reversed} along one of the edges being twisted. }
\end{figure}

Beyond the sedenions, twists no longer always take ZD edges to ZD
edges.  Type I's always twist to Type I's; but Type II's, first
emerging in the pathions, either twist to other Type II's, or to
box-kite-like structures \textit{none} of whose edges act as ZD
pathways.

Per the Roundabout Theorem of [2], box-kites are ``all or nothing''
structures:  all edges support ZD-currents, or none do. These latter
\textit{hidden box-kites} (HBKs, or ``residents of Hoboken'') were
the sources of the off-diagonal empty cells in the $2^{N - 1}- 2$
cells-per-edge square ETs for fixed-{\bf S} $2^{N}$-ions, studied in
[2],[3], and presented in color-coded spreadsheet displays in our
NKS 2006 Powerpoint slideshow [9]. These showed what [2] and [3]
proved: that, as $N$ grew indefinitely large for fixed {\bf S}, such
tables' empty space approached a fractal limit.

For $N = 4$, each of the 7 ETs is a 6 x 6 table, one label per each
possible L-index excluding {\bf S}; for $N = 5$, {\bf S} takes all
integer values less than ${\bf G} = 16$, with edge-length in each ET
being 14 (the number of indices $< {\bf G}$, with {\bf S} excluded).
Consider $N = 4$, and ignore the $2 \cdot 6$ cells along long
diagonals:  these are tautologically empty, since ZDs in the same
assessor do not mutually zero-divide, nor do those of assessors
which are strut-opposites.  24 filled cells remain: two for each
edge, hence one for each distinct ZD-pairing defined on it. This
shows the ET is fundamentally a \textit{multiplication table}, with
only L-indices indicated on the row and column headers, in
nested-parentheses order (i.e., the leftmost assessor label A is
strut-opposite to the rightmost label F, then B to E, and so on by
mirror symmetry). This is because U-indices are forced, hence can be
ignored, for given {\bf S} and $N$.

For any assessor $(M, m)$ and its strut opposite $(M_{opp},
m_{opp})$, it's easy to see that $m \veebar ({\bf G + S}) = {\bf G}
+ m_{opp} = M$.  Twist products along an edge are hence linked with
a box-kite whose {\bf S} is the L-index of the assessor which is the
strut opposite of the third assessor in the given edge's sail. Both
$(A + b)$ and $(B + a)$ then have ${\bf S} = d$, since $A$ and $B$
are in a sail with $C$, whose strut-opposite is $D$. And by the
third production rule, we know the edge opposite that joining $A$
and $B$ also has its sail completed by $C$:  that is, the square
formed by $(A, B, E, F)$ and orthogonal to the strut $(C, D)$ will
have 4 of the 6 assessors defining the box-kite with ${\bf S} = d$
residing along one set of parallel sides, while 4 of the 6 defining
the ${\bf S} = c$ box-kite will reside along the other parallels in
the same square.  (Corollary: for any box-kite, each L-index is also
the {\bf S} of another box-kite reached by twisting.)

With 3 such \textit{catamarans} per box-kite, each with edges whose
sails are completed by assessors of a different strut, all 7
sedenion box-kites can be seen as collected on the frame of just
one. The missing pairs of assessors are derived by twisting the $(
{\bf S}, {\bf G + S}) \equiv ( {\bf S}, {\bf X})$ pair, imagined in
the center, with each of the legitimate assessors, yielding
assessor-pairs defined along each catamaran's mast and keel
(strut-halves (a, A) to ({\bf S, X}), then ({\bf S, X}) to (f, F),
in that order, in Figure 1.) Such a 7-in-1 representation we call a
\textit{brocade}.

In the table, the singleton sedenion brocade shows all possible
L-indices as column heads, U-indices as row labels, and a long
diagonal of empty cells signifying the ({\bf S, X}) non-assessor
pairs.  Each cell gives {\bf S} and vertex letter for all 42
assessors specified by U- and L- indices. The zigzag for ${\bf S} =
1$, say, is $(3, 10); (6, 15); (5, 12)$, with twists $(b, A) = (6,
10)$ and $(a, B) = (3, 15)$ yielding assessors E and C of the ${\bf
S} = 4$ box-kite. For $N > 4$, seven pairs of row and column labels
still fix one brocade, but indices will no longer be consecutive,
and cellular information will need to indicate which of the numerous
box-kites is being twisted to among those of all types with the same
{\bf S} -- a number equal to the trip count in the $2^{N - 2}$-ions
given earlier, which surprising result was derived as a corollary of
the Roundabout Theorem in [2].

\medskip

\begin{center}
Table 1:  The Sedenion Brocade \linebreak

\fbox{
\begin{tabular}{ c | c | c | c | c | c | c | c }
  $\;$  & ${\bf 1}$   & ${\bf 2}$  & ${\bf 3}$  & ${\bf 4}$  & ${\bf 5}$  & ${\bf 6}$  & ${\bf 7}$  \\ \hline

  \textbf{09}    & $\;$  & $3A$  & $2F$ & $5B$  & $4F$  & $7F$ & $6C$ \\ \hline
  \textbf{10}    & $3F$  & $\;$  & $1A$ & $6B$  & $7C$  & $4E$ & $5F$ \\ \hline
  \textbf{11}    & $2A$  & $1F$  & $\;$ & $7B$  & $6F$  & $5C$ & $4D$ \\ \hline
  \textbf{12}    & $5E$  & $6E$  & $7E$ & $\;$  & $1C$  & $2C$ & $3C$ \\ \hline
  \textbf{13}    & $4A$  & $7D$  & $6A$ & $1D$  & $\;$  & $3E$ & $2B$ \\ \hline
  \textbf{14}    & $7A$  & $4B$  & $5D$ & $2D$  & $3B$  & $\;$ & $1E$ \\ \hline
  \textbf{15}    & $6D$  & $5A$  & $4C$ & $3D$  & $2E$  & $1B$ & $\;$ \\
\end{tabular}
 }
\end{center}

\medskip

Our prior work showed that an ET's empty cells -- emerging in any
and all ET's for $N > 4$, and ${\bf S} > 8$ and not a power of $2$
-- mapped to pixels in a planar fractal. Here, catamaran twisting
will let us see how these emptinesses have their own subtle
structure, coming in quartets of two distinct classes, exactly akin
to zigzags and trefoils among sails.  Moreover, any box-kite is
uniquely linked to four HBKs, with each such quartet or
\textit{spandrel} housing its own ZD-free copy of the octonions –
hence the basis for the recursive CDP spawning of independently
generated $2^{N}$-ion index sets, or \textit{Context-Definition
Platforms} (whose acronymn's meaning is, of course, itself
context-dependent). \pagebreak

\section{ Box-Kite ``Explosions'' in 32D:  Two Types, Triptych Triples, Four-Fold Spandrels }

Historically, a famous proof from the late 1890's by Adolf Hurwitz
[10] dissuaded researchers from investigating any $2^{N}$-ions
beyond the sedenions:  once Hurwitz showed that they, and all higher
hypercomplex numbers, unavoidably contained ZDs, the entire line of
study was deemed pathological -- hence, our calling those in 32-D
(the smallest-$N$ $2^{N}$-ions to not have a name) the
\textit{pathions}. But, as with their contemporary ``monstrosities''
of analysis, whose taming by Benoit Mandelbrot led to fractals, the
pathions in fact mark the beginning of a new agenda, at least as
much as they signal the demise an older one.  The work just prior to
this paper shows that the connection to Mandelbrot's discoveries is
not just by analogy:  as a side-effect of what we might think of as
carry-bit overflow, ETs in high-N $2^{N}$-ions, beginning with the
pathions, have surprising patterns of empty cells when {\bf S} is
not a power of 2, and its binary representation contains one or more
bits to the left of the 4-bit.

In the pathions, 15 L-indices $<$ {\bf G}, hence candidate {\bf S}
values, times 7 (the octonion trip count) per ET, means 105
box-kites. Seven are the equivalent of those found in the sedenions,
but for the zero-padding of {\bf G} (via left-shifting its singleton
bit), and hence of {\bf X}: all L-trips are identical, but U-indices
at each assessor are augmented by the difference of the old and new
{\bf G} values, or 8. The 7 box-kites for ${\bf S} = 8$ (the
sedenions' {\bf G} ) are Type I, but are special in other regards.
First, the Z-trip of each is the same as one of the sedenions';
hence, these seven Rule 0 trips, once {\bf S} is downshifted to its
sedenion twin's value, can map directly to one of the zero-padded
box-kites.  Similarly, each \textit{strut} is a Rule 1 trip, serving
as the $(a, d, e)$ L-trip of a pathion box-kite, with the same
downshifted {\bf S}.

Finally, the 3 trefoil L-trips are just Rule 2 transforms of the
Z-trip (since this ${\bf S} = 8$ acts on it as a minimal {\bf G}).
Z-trips in their own right, they also produce box-kites with
downshifted {\bf S} values -- of the new Type II.  We thus have at
least $7 \cdot 3 = 21$ of these in the pathions. In fact, we have
\textit{only} these 21, derived from trefoil L-trips of ${\bf S} =
8$ box-kites; hence, the add-and-switch logic of Rule 2 should be
central to their new typology.  As is the case:  exactly 2 of the 3
struts in a Type II have their orientations reversed, as mentioned
above.  Each Z-trip index of a Type I gives its strut-opposite
assessor's L-index when multiplied on the right by {\bf S}, but 2 of
the 3 Type II zigzag's L-units form CPO struts when multiplied by
{\bf S} \textit{on the left}. \pagebreak

As shown in Figures 3 and 4, we can visualize all this by adapting
the commonplace Fano plane rendering of our XOR-based octonion
labeling scheme to different ends, a.k.a. the PSL(2,7) Triangle --
for "projective special linear group of 7 lines in the plane,''
which cross in 7 places.  This simplest nontrivial finite projective
geometry has each line projectively equivalent to a circle -- which,
adapting standard convention, is how only the Rule 0 Z-trip is
drawn.  The 3 lines through the central node join angles to
midpoints, making Rule 1 trips when the label in the center is a
power of 2. The 3 sides then become the Rule 2 trips, in the manner
just discussed:  the center is the \textit{sedenion} {\bf G},
converted to a pathion {\bf S}.

The left of Figure 3 can be read as displaying the L-trips of the
sedenion box-kite with ${\bf S} = 4$; one inflates the diagram by
assuming the attaching of U-indices, by the $L \veebar {\bf X}$
rule, to get a full box-kite, each side now turned into a bonafide
trefoil sail, and the Rule 0 L-trip converted to a zigzag. The
diagram at the right abstracts this via assessor L-index lowercase
coding conventions.  The approach just sketched works for
shorthanding box-kite structures for any $2^{N}$-ions, $N \geq 4$.

Note the CPO flow along all lines:  the triangle's perimeter is
naturally traversed clockwise, as is the central Rule 0 circle,
while strut-flows move from midpoints, through {\bf S}, to the
angles. In the next two diagrams, the sedenion Z-trip for ${\bf S} =
1$ doubles as the pathion Z-trip for one of the 7 ${\bf S} = 8$
box-kites; then, one of its Rule 2 sides is inflated on the right,
to yield a pathion ${\bf S} = 1$ box-kite of Type II.  Note its flow
reversals along the Z-trip's $b$- and $c$- based struts.

\begin{figure}[H]
\centering
\includegraphics[width = 4in, height = !] {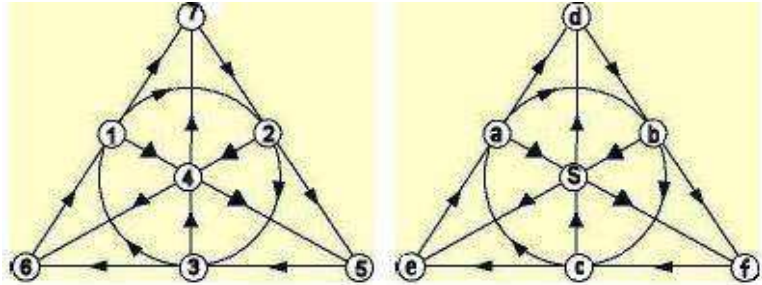}
\caption{ The 7-point, 7-line finite projective group, a.k.a the
Fano Plane, hosts the labels for octonion units, and shows their
triplets' orientations.  The same layout can be used to shorthand
box-kite structures:  the zigzag and trefoil L-trips sit at ${\bf
(a,b,c)}$, and $({\bf a},d,e)$; $(d,{\bf b},f)$; $(e, f,{\bf c})$.
The strut constant {\bf S}, meanwhile, sits in the middle.}
\end{figure}

\begin{figure}[H]
\centering
\includegraphics[width = 4in, height = !] {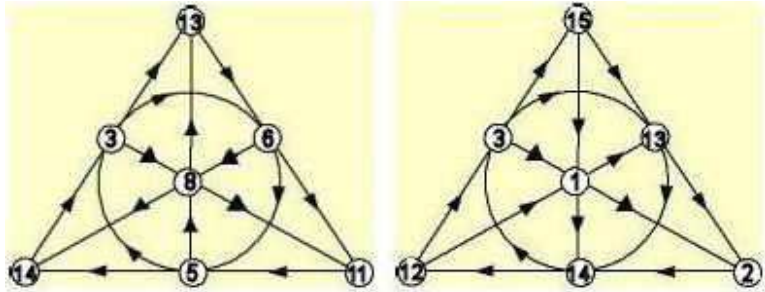}
\caption{Pathion box-kites:  Left, a normal (Type I) with ${\bf S} =
8$, and Rule 0 Z-trip $(3,6,5)$ at $(a,b,c)$ -- itself the Z-trip
for the ${\bf S = 1}$ sedenion box-kite. Right, a Type II with ${\bf
S} = 1$, with Z-trip the Rule 2 left side of the ${\bf S} = 8$ Type
1.}
\end{figure}

In Theorem 7 of [1], the parallel flows around the triangle's
perimeter and central circle provided the implicit basis for proving
the PSL(2,7) in question was a Type I box-kite.  What we now call
the sedenion \textit{brocade} compactly expresses the fact that,
provided the node-to-node connections and flow patterns aren't
changed, any node can be moved into the center to act as the strut
constant, with the only substantive side-effect being the
broad-based swapping of U-indices associated with each node.

Direct hand calculation makes it clear that Theorem 7 still holds
for a Type II box-kite, as the flows remain parallel around the
Triangle's perimeter and inscribed circle.  A Type I twists to
another Type I.  A Type II, though, only twists to another Type II
when the single strut with proper orientation (all of whose
L-indices, in the pathions, are octonions) has one of its nodes
swapped into the center (or, equivalently, provides the strut for
the catamaran being twisted). For all other twistings, ${\bf S} >
8$, and we have HBK's -- tantamount to saying (albeit not in an
obvious way) that the perimeter and circular flows no longer stay
parallel.

In Theorem 15 of [3], we proved that two L- and U- unit pairings
which mutually zero-divide (share an edge as assessors on a
box-kite) no longer do so once {\bf S} is augmented by a new high
bit:  this was the general case inspired by the empirical
for-instances provided by the pathions' ETs (for $8 < {\bf S} <
16$). Here, only 3 of the 7 box-kites for any such {\bf S} prove to
be Type I; the remaining 4 reside in Hoboken.

Using the theorem just cited, each such {\bf S} is just that of a
sedenion box-kite with the minimal new high-bit appended to it. This
makes strut-opposite L-indices, whose XOR is {\bf S}, have their
difference augmented by 8 -- which means the larger are U-indices of
the sedenion case.  We can then take each assessor in a sedenion
box-kite and treat it as a pathion pair of L-index strut-opposites,
effectively exploding one assessor into two.

This implies each sail can be inflated into its own box-kite,
sharing one strut with each box-kite built, by the same logic, from
each other sedenion sail.  And, as the theorem will apply similarly
to each and every sedenion edge-current, and hence all 4 of the
L-trips, we can say each of the sedenion L-trips does service as a
Z-trip for an HBK.  We call such quartets of Hoboken residents
\textit{spandrels}, after a term made famous by evolutionary
theorists Steven J. Gould and Richard Lewontin [11].

The deep appropriateness of this term will become apparent in the
final remarks of this paper, when we consider the epistemological
issues it was coined to address.  But a superficial aptness is
easily grasped.  Consider one of the secondary meanings of the term
(which Gould and Lewontin didn't have in mind):  among philatelists,
the four curved wedges between perforated border and inner oval
containing, say, a president's face, comprise a postage stamp's
\textit{spandrel}.  Pinch diagonally opposite corners of such a
stamp together, so that two meet above, and two below, the center of
the stamp proper.  The spandrel's wedges become sails in a Box-Kite
(the kind of corner-to-corner mapping of flows on a plane, from
which the \textit{projective plane} derives).

Each of the 7 sedenion box-kites explodes into one pathion spandrel
-- 28 HBK's in all.  Simple arithmetic shows how this count
dovetails with what was said above about Type II's: one can only
twist to two other Type II's from a given one, their strut constants
forming an octonion trip with that we are starting with. All 4 other
twists take one to Hoboken, where each box-kite has one all-octonion
L-trip inherited from its sedenion box-kite of origin. Its own {\bf
S} being larger than the prior {\bf G}, it can be twisted to 3
different Type II's, hence 3 other HBK's. Ergo, there will be 3 Type
II's for every 4 HBK's -- or 21 for the 28 HBK's in the pathions, as
already calculated.

But there will also be 3 Type I's, each of whose BE strut comprises
the assessors whose L-indices are the former {\bf G} and {\bf S} of
the sedenion box-kite they exploded from.  Further, each former
strut now has its vent and zigzag L- (and U-) indices appearing at
$a$ and $d$ (and $f$ and $c$) respectively (forming trefoil L-trips
thereby with the old {\bf S} and {\bf G} at $e$ and $b$), in one of
the 3 new pathion box-kites. These trios are the \textit{sand
mandalas} first reported on (and graphically rendered) in [12],
which we generalize to the general $2^{N}$-ion case by redubbing
them lowest-N examples of \textit{triptychs}.

In the general case, however, while 3 box-kites are exploded from
each Type I we start with, they are not unique in derivation.  Each
corresponds to a $2^{N}$-ion strut that's been inflated into a $2^{N
+ 1}$-ion box-kite. But the $2^{N + 1}$-ions have $2^{N} - 2$
distinct assessors (hence a strut-count half that number) shared
among all same-{\bf S} box-kites (including Hoboken residents).
Hence, we do indeed get 3 for each pathion ET with ${\bf S} > 8$.
But, for any zero-padded sedenion (hence, pathion) Box-Kite
(meaning, {\bf G} is left-shifted to be 16 not 8), although we might
start with the ${\bf S} = 1$ case, it now houses \textit{7} distinct
box-kites, not 1, for all of which {\bf X} is 17 instead of 9.  The
results of explosion thereby reside in the chingons, not pathions,
where the ET for ${\bf S} = 17$ has 15 struts shared among seven
distinct box-kites (one of which struts, with L-indices \textit{b} =
16 and \textit{e} = 1, is shared by all seven box-kites), not (as in
the pathion case of the septet of "sand mandalas" with $8 < {\bf S}
< 16$) seven struts shared among three box-kites (with, e.g.,  the
BE strut \textit{b} = 8 and \textit{e} = 1 held in common when ${\bf
S} = 9$). As with the pathion ${\bf S} = 9$ case, though, these are
the ET's \textit{only} non-empty box-kites: the remaining 4 x 7 = 28
are all contained in pathion-generated spandrels.

The Greek etymology of ``triptych'' indicates 3 (\textit{tri-})
plates or panels (\textit{ptyche}).  The `tri' indicates the count
of BE-sharing distinct box-kites in the pathion case only; more
generally, the `panels' are box-kites associated with the
(pre-explosion ET's) distinct strut triplets instead -- and the
count of these `trips' in a triptych is typically much higher than
3. In this sense, they are akin to what a Java programmer might call
``static variables'':  unlike the spandrel quartets, their
generation is tied to the ET ``class'' of origin, rather than to a
particular source box-kite.

Described in this manner, triptychs may seem more concocted than
natural.  This is not so when viewed from a purely bit-twiddling
vantage:  when, as shown in Figure 5, their ETs are examined, the
flip-book sequence generated by integer increments of S between 8
and 16 shows animation logic. Four lines just off the picture frame,
spanning the long diagonals' empty corners, form the 12-cell-long
sides of a square including the corners, hence taking up the maximum
$14 \cdot 14$ size that a pathion ET allows. (Similar descriptions
obtain for the $30 \cdot 30$-sized chingon flip-books for $16 < {\bf
S} \leq 24$.)

As {\bf S} grows, these orthogonal pairs of parallels move one cell
in from the perimeter with each increment, until, when ${\bf S} =
15$, they form two-ply cross-hairs partitioning the ET into
quarters.  The remaining 24 filled-in cells form 6-cell-long
diagonal spans, connecting the cross's vertical and horizontal ends.

This abstract cartoon or \textit{flip-book} is drawn by a simple
formula, the gist of Theorem 14 in [3]:  using the vertical pipe for
logical OR, and shorthanding the {\bf G} of the $2^{N-1}$-ions as
$g$,

\begin{center}
${\bf R\; |\; C |\; P} = g\; |\; {\bf S} \; mod \; g$
\end{center}

\begin{figure}[H]
\centering
\includegraphics[width = 4in, height = !] {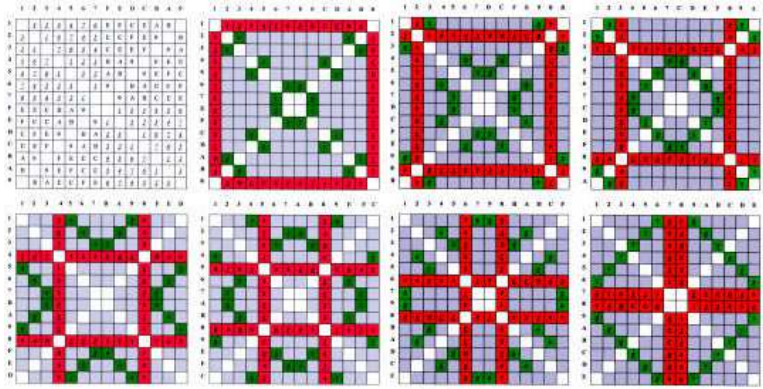}
\caption{8 pathion ET's with ascending {\bf S} values, from the
maximally full ${\bf S} = 8$ to $15$ on the second row's far right.
For all ${\bf S} > 8$, the 24 cells filled with the sedenion {\bf G}
and {\bf S} are painted darker than the 48 forming the parallels.
Two shades of grey highlight the $4 \cdot 24$ long- and off-
diagonal unfilled cells of the spandrel's HBK's.}
\end{figure}

Only if row label {\bf R}, column label {\bf C}, or their XOR
product {\bf P}, equal $g$ or $s$ (the pre-explosion {\bf S}), will
the cell be filled (and assessors with L-indices {\bf R} and {\bf C}
mutually zero-divide). Via Recipe Theory [3], this formula can be
generalized, by a simple analysis of {\bf S}'s bitstring. For any
${\bf S > 8}$ not a power of 2,  the ETs' empty spaces for each
successive $N$ approach a fractal, overlaying each other's values.
Row and column \textit{labels} of the $2^{N-1}$ ET become actual
cell \textit{values} of the $2^N$, with the same values filling in
the label-lines' empty parallels in reverse (strut-opposite) order,
in a never-ending \textit{balloon ride} sequence (see Figure 6) of
nested \textit{skyboxes}.

\begin{figure} [H]
\includegraphics[width = 1.\textwidth] {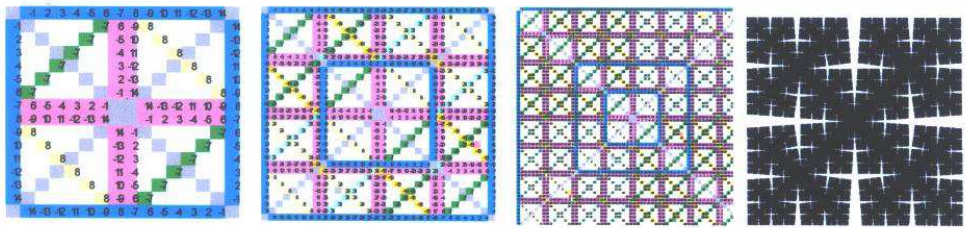}
\caption{ETs for S=15, N=5,6,7 (nested skyboxes bordered in darker
shading)$\cdots$ and fractal limit:  the Ces\`aro double sweep [13,
p. 65].}
\end{figure}

In this sense, Recipe Theory is a pure Wolfram-style number theory,
focused on the binary representations of integers rather than their
quantities, hence according special status to the placeholding power
of \textit{singleton bits} ({\bf G} values) -- as opposed to
\textit{traditional} number theory, which concerns itself, above all
else, with \textit{size} -- and hence, with \textit{primes}.

Complementary to Theorem 15 of [3], just cited above, is the Theorem
16 which immediately follows:  while a box-kite's edges are turned
off by augmenting its {\bf S} with a new leftmost bit (and
necessarily left-shifting its {\bf G}-bit if this new {\bf S}
exceeds it), performing a second such augmenting results in a
box-kite which is once again turned on.

In addition to the \textit{$(s,g)$ modularity} first seen in our
sand-mandala formula, we thus have a process of \textit{hide/fill
involution}:  repeated, it produces spandrels from proper (Type I or
Type II) box-kites; quartets of higher-N proper box-kites from each
such HBK; quartets of higher-N spandrels from each of
\textit{these}; and so on, ad infinitum. This is a result
sufficiently powerful as to call for a proof.

For the HBK deriving from a Type I's zigzag, the trips along all 3
struts are reversed:  if $(z, {\bf S}, v)$ is CPO, then Rule 2 says
replacing {\bf S} with {\bf X} by exploding the Type I's assessors
will reverse orientations.  Similarly, the flows along the edges
will also reverse, since each edge's two terminal nodes are
reversed.  So if an appropriate power of 2, $g$, be added to the
central node, this same $g$ must be added to the 3 nodes at the Fano
plane's angles as well, to keep all lines trips. Only one among the
seven Fano lines will have no, instead of two, nodes with $g$
appended to its indices:  the Z-trip itself.  Hence, if the Z-trip
flows clockwise, perimeter tracings now run counterclockwise.

Such a fourfold $g$ insertion is equivalent to exploding a
$2^{N}$-ion box-kite's zigzag, since $g$ appended to the strut
opposite $v$ of any Z-trip L-index $z$ yields its pre-explosion
U-index partner:  for arbitrary Type I zigzag assessor (z,Z), $z
\cdot {\bf S} = v$, but $v \cdot ({\bf S} + g) = v \cdot {\bf X} =
Z$.  The resulting Fano's edge trips, in terms of pre-explosion
indices, now read (B, a, C); (C, b, A); (A, c, B) -- all with
orientations reversed.  For recall a crucial fact from earlier work,
which we termed ``trip-sync'':  in the zigzag only, the 4 quaternion
copies $(a, b, c); (a, B, C); (A, b, C); (A, B, c)$ -- the L-trip
and its 3 allied \textit{U-trips} -- all flow similarly, so that one
can effectively allow ``slippage'' between high- and low- index
units at any of the zigzag's assessors and not notice any
difference. Repeating this process with fourfold insertion of ${\bf
G} = 2 \cdot g$ re-reverses all six flipped lines, again leaves the
Z-trip unchanged -- again giving a Type I box-kite.

For the 3 HBK's derived from making a Z-trip out of a Type I's
trefoil L-trip, the 3 trefoil-derived HBK's in any Type I box-kite's
spandrel will have \textit{just one} strut-and-edge pairing reversed
-- and a different one in each case. This is the complement to the
zigzag ``slippage'' effect:  among the 4 quaternion copies
associated, say, with $(a, d, e)$, only $(a, D, E)$ will share its
orientation:  the other 2 U-trips, whose singleton L-index is not
shared with the zigzag, have CPO forms $(A, E, d)$ and $(A, e, D)$.
Each of the trefoil HBK's, then, has flow structure like that of a
Type I save for a ``T'': the reversed strut flows from an angle to
the midpoint of the likewise flow-flipped side-trip.  Doing the
fourfold $g$-appending to the nodes \textit{not} included in the
only reversed line (containing the only L-index from the
pre-exploded box-kite's zigzag) will also result, as direct symbolic
calculation shows, in a Type I box-kite.

All told, then, if none or two of the struts are reversed, we have
proper box-kites, of Types I and II respectively; if all or one, we
have zigzag or trefoil HBK's in that order for Type I spandrels.
Type II spandrels appear, at first, to have their own distinct Fano
Plane flow-patterns.  However, by swapping sides with zigzags, the
HBK reversal patterns take on very different appearances -- and the
full set of 16 sail-based presentations are the same in both
spandrels, but in different orders.  The display with central circle
housing the default $(a,b,c)$ line, and that with $(a,d,e)$, reverse
places in Type I and Type II spandrels.  This shuffling is crucial:
because of it, \textit{both} types will be found to facilitate
\textit{cowbirding}, a process of paramount interest.

Successive explosions of the sedenions' ${\bf S} = 1$ box-kite take
us to the ${\bf S} = 9$ and ${\bf S} = 25$ ET's of the pathions and
chingons respectively.  Such instantiating is readily generalized,
since each \textit{spectral band} of 8 consecutive {\bf S} values
(powers of 2 excluded) obeys the same hide/fill logic: like the 7
sand-mandala ET's in the pathions, there is an animation-like
impetus connecting each to each, all with the same counts of proper
and hidden box-kites.  For those who like to read the libretto at
the opera, the graphics corresponding to the cases just mentioned,
in the order just given, are on Slides 16, 25, and 48 of [9].

Our sedenion starter kit has Z-trip $(3,6,5)$.  For this and its
L-trips, we get the 4 HBK's, all with ${\bf S} = 9$, ${\bf X} = 25$,
${\bf G} = 16$, as shown below.  (The original assessor L- and U-
indices, all with ${\bf S} = 1$, ${\bf X} = 9$, ${\bf G} = 8$, are
shown on the first and third lines; as pathion L-indices, their
corresponding U-indices appear in the second and fourth.)

The zigzags of the pathion HBK's share a special feature: treat all
6 L- and U- indices, plus {\bf X} and the Real unit, as a set.  Not
only are the HBK's edges bereft of ZD currents; this 8-D ensemble
shows no ZD currents anywhere within it, no matter how you twist it.
No two pairings of one element each from $(a,b,c)$ and $(A,B,C,X)$
will mutually zero-divide. We have, then, a pure octonion copy.
\pagebreak

\begin{center}
Table 2:  Sedenion Sails and their Pathion Explosions \linebreak
$({\bf S, X, G} = 1, 9, 8) \Rightarrow ({\bf S, X, G} = 9, 15,
   16)$ \linebreak

\fbox{
\begin{tabular}{ c | c | c | c }

  $(a,b,c)$    &    $(a,d,e)$     &    $(d,b,f)$    &    $(e,f,c)$        \\ \hline

  $(03,06,05)$ &    $(03,04,07)$  & $(04,06,02)$    &    $(07,02,05)$     \\
  $(26,31,28)$ &    $(26,29,30)$  & $(29,31,27)$    &    $(30,27,28)$     \\ \hline

  $(10,15,12)$ &    $(10,13,14)$  & $(13,15,11)$    &    $(14,11,12)$     \\
  $(19,22,21)$ &    $(19,20,23)$  & $(20,22,18)$    &    $(23,18,21)$     \\

\end{tabular}
 }
\end{center}

\medskip

Clearly, the four Q-copies involving L- and U- units along the
zigzag's edges are ZD-free, as are the 3 trips involving the
spandrel's {\bf X} with the U- and L- indices of each zigzag
assessor in turn.  Such an ensemble, whose indices can be relabeled
to be identified with standard octonions, we call an \textit{egg};
its habitat, a \textit{cowbird's nest}.  Such nidi exist in all
spandrels, one per each of the four HBK zigzags, in their zigzag
(ADE) sails for Type I (Type II). Implied in this are three claims,
comprising the \dots

\medskip

\noindent {\small Cowbird Theorem.}  Each HBK in a spandrel contains
a ZD-free octonion egg $O$, consisting of the Reals, the L-indexed
units of one of the source box-kite's sails, their U-index partners
in the HBK exploded from it, plus another imaginary indexed by the
spandrel's {\bf X}. $O$-containment is \textit{universal} among HBKs
found in either Type I or Type II spandrels, and is also
\textit{habitat-specific}: for Type I, an $O$ is always and only to
be found in an HBK's zigzag sail, whose L-index set is identically
that of the source box-kite's sail from which the HBK was exploded.
Type II spandrels, meanwhile, harbor eggs in their HBK's ADE sails,
in patterns that are \textit{flowmorphic} -- that is, have
identically connected and oriented lines in their Fano presentations
-- to the Type I eggs' nests.

\medskip

We introduce and exemplify seven lemmas that pave the way to our
proof. These highlight and customize familiar aspects of PSL(2,7),
the 168-element simple group governing manipulation of the Fano
diagram's labeling.  (This is not to be confused with the 480
distinct labeling schemes for octonions mentioned earlier:  only a
small subset of these entail the XOR relations defining triplets
upon which our apparatus depends.)  We also want to indicate certain
correspondences between the fourfold {\bf G} \textit{loading} -- at
the corners, and in the center -- which effects the explosion
process taking one from $2^{N}$ to $2^{N+1}$ levels, and the
same-level fourfold \textit{label-exchange} process which
underwrites the basic PSL(2,7) symmetries.

\medskip

\noindent \textit{Lemma 1}:  Explosion fixes one of the 7 Fano-plane
lines (the zigzag circle, in our arguments' context), with the other
6 having their orientations reversed:  all 6 contain exactly 2 of
the 4 high-bit-augmented nodes, so Rule 2 applies to each. By Fano
symmetry, the same holds true whatever four nodes are sites for {\bf
G} loading, provided no 3 reside in one line (occurring only if
exactly one line has no nodes selected).

\medskip

\noindent \textit{Lemma 2}:  If the fixed line be one of the 3
sides, the G-loading nodes will form a kite containing {\bf S} and
the angle opposite the fixed side, and the 2 midpoint nodes joined
to both. If a strut, G-loading sites are the node-pairs
perpendicular to it, residing in midpoints and vertices. But all 3
cases are projectively equivalent, and any 4-node loading can be
uniquely specified writing \textit{PL($m,n,p$)}, triplet $(m,n,p)$
the fixed line, each line choice uniquely corresponding to the
explosion of one of the 7 (possibly hidden) box-kites in a single
brocade.

\medskip

\noindent \textit{Lemma 3}:  Consider two distinct lines $l_{1}$ and
$l_{2}$ among a Fano presentation's 7. If $N(P)$ lines $\in$ P, the
set whose orientations are \textit{preserved} (oriented identically
to those placed in the same positions in a standard Type I box-kite,
per Figure 3), and $N(R) = (7 - N(P))$ lines $\in$ R, where
orientations are \textit{reversed}, perform PL($l_{1}$). For $l_{1}$
only, orientation stays unchanged.  If $l_{1}$ $\in$ P, then after
the performance, $N(R) = (N(P) - 1)$; if $\in$ R, then $N(R)= (N(P)
+ 1)$. If we start with standard presentations of either Type I
($N(R) = 0$), Type II ($N(R) = 2$), Type I's explosion ($N(R) = 6$)
or Type II's ($N(R) = 6$ if the fixed line $l_{1}$ which produced it
begins in R, else $N(R) = 4$), then any followup performance of
$PL(l_{2})$ will have $N(R)$ an even number.

\medskip

\noindent \textit{Lemma 4}:  The four-point sites for G-loading are
the necessary and sufficient bases for \textit{double exchange} (DX)
of labels, the minimal non-trivial label exchanges possible which
preserve the mutual connections among the 7 lines. For a fixed line
$(m,n,p)$, the expression \textit{DX($m,n,p$; $Op$)}, $Op$ either
$H,$ $V$, $D$, or the identity operator $I$, defines a Klein Group.
Taking the fixed line to be the vertical strut ($d, {\bf S}, c$) on
the standard Type I layout for ease of visualizing, the horizontal
DX operator $H$ (instantiated here via $a \rightleftharpoons b$, $e
\rightleftharpoons f$) produces a dihedral flip of the Fano triangle
along $(d, {\bf S}, c)$.  The vertical DX $V$ (realized here as $e
\rightleftharpoons a$, $f \rightleftharpoons b$) makes the trefoil
$(e,c,f)$ and zigzag $(a,b,c)$ trade places while preserving their
(and the vertical strut's) orientations. The composite $D$ = $H
\cdot V = V \cdot H$ is commutative, and exchanges diagonally
opposite members of the four-point site.  All three operators being
involutions, their Klein Group can be written by supplementing the
$D$ composition just given with the relations $H^{2} = V^{2} = D^{2}
= I.$ (As $168 = 4 \cdot 7 \cdot 6$, with 7 lines of triplets each
permutable in 3! = 6 ways, the structure of the Fano plane's group
qua \textit{brocade presentation} is hereby given.)

\medskip

\noindent \textit{Lemma 5}:  Define counts $N(R)$, $N(P)$ of lines
in sets R and P as in Lemma 3.  Assume a given Fano presentation is
derived from standard Type I or Type II presentations by PL
operations, as in Lemma 4.  Then $N(R)$ is even. Consider the before
and after counts for the dihedral flip along the vertical strut,
DX($a{\bf S}f; V$).  Orientations of all 3 struts are unaffected;
hence, if their contribution to $N(R)$ before DX is odd, it will be
so after DX, and ditto for DX even.  Orientations for all 3 sides
plus the center, meanwhile, are all reversed:  hence, if their
contribution to $N(R)$ before DX is $k,$ it will be $(4 - k)$ after
DX, hence its odd or even status will remain unchanged.  If we
start, then,  with Type I or Type II or any HBK derived from either
by multiple applications of PL, $N(R)$ will be even before DX, and
remain so, since odd plus odd and even plus even are both even.  As
all DX operations, per Lemma 4, are equivalent in all senses
relevant, all DX operations preserve evenness of the reversed-flow
count.

\medskip

\noindent \textit{Lemma 6}:  Any standard-form Type I can be
converted into any of 3 different standard-form Type II's by
performing, in either order, a DX and a PL on a strut $m{\bf S}n$.
For the zigzag $(3,6,5)$ of the ${\bf S} = 1$, $N = 4$ (hence, Type
I) box-kite, DX($a{\bf S}f; V)$ and PL($a{\bf S}f)$, performed in
either order, gives two easy ways to generate the Type II shown on
the right of Figure 4.

\medskip

\noindent \textit{Lemma 7}:  By Lemmas 3 and 5, any combination of
PL and DX operations (mutually commuting where
orientation-preservation is concerned), will always reverse an even
number of lines; further, the number of reversed lines with respect
to the standard Type I presentation will likewise always be even
after performing any number of such operations. Ergo,
standard-presentation Type III and Type IV box-kites (Type I with 1
or 3 struts reversed respectively) can never be derived from Type I
and Type II types.  But G-loading, and catamaran twisting within a
brocade, are the bases for all possible creation of standard
box-kites.  Hence Types III and IV are impossible.

\medskip

Now let's tackle the cowbirding theorem, give reasons for so naming
it, then spell out our claims for its potential to underwrite models
of semantic networks.

\medskip

\noindent \textit{Proof.} We consider only Type I for most of our
argument, then find a simple way to carry over everything shown
therein to the Type II situation.  Label a standard-form Type I
box-kite in the $2^{N}$-ions in conformity with that shown on the
right of Figure 3: strut constant {\bf S} appears in the center, the
zigzag's $a,b,c$ are arrayed about it at 10, 2, and 6 o'clock
respectively, with their strut opposite nodes $f,e,d$ appearing in
the lower right and left corners and apex, in that order.  Explode
by G-loading the corners and center with $g$ (the pre-explosion {\bf
G}) -- synonymous with saying perform the operation PL($a,b,c$). The
result is the $(a,b,c)$-based HBK, with {\bf X} now equal to the
pre-exploded ${\bf S = s}$, and new ${\bf S} = s + g$, and new ${\bf
G} = 2g = 2^{N}$. To prevent ambiguity when pre- and post- explosion
nodes are referenced, the subscript `H' designates nodes in the HBK,
not the source box-kite.

For arbitrary zigzag assessor $(Z,z)$ and its strut-opposite vent
assessor $(V,v)$ in the HBK, U-indices $Z$ = ${\bf G} + v$; but
PL($a,b,c$) having loaded each pre-exploded $v$ with $g$, this means
$(z,Z) = (z, {\bf G} + g + v)$, and $(v,V) = (v, {\bf G} + g + z)$.
Our initial claim: in this HBK, the 6 zigzag units, {\bf X} $( =
{\bf G} + g + s)$, and the Real unit (index $= 0$) are isomorphic to
the standard octonions, hence ZD-free.  Call this octet \textit{egg}
$O$ in the $(a,b,c)$ cowbird's nest, associated with Type I
box-kite, having ${\bf S} = s$ and zigzag L-index set $(a,b,c)$ in
the $2^{N}$-ions, shorthanded $(a,b,c; s; N)$.

The multiplication table of $O$ is clearly isomorphic to that for
the usual octonions: for each $(z,Z)$ choice, $z \cdot {\bf X} = Z,$
and HBK trips $(a,b,c)$, $(a,B_{H},C_{H})$, $(A_{H},b,C_{H})$, and
$(A_{H},B_{H},c)$ provide the remaining four quaternion copies that
complete the septet found in standard octonions.  And, as this is an
HBK, we already know that the products represented by the HBK's
zigzag edges are \textit{not} zero.  We now take a ``voyage by
catamaran,'' showing, by twist-product logic, how the other L- and
U- index pairings contained in $O$ \textit{also} make no zeros.

For each pairing among zigzag L-indices, the Fano triangle
presentation makes it clear that the arrow sweeping out $120^{o}$
between them forms an arc beneath the corner L-index which
represents the strut constant of the box-kite the assessors
containing said L-indices twist to.  For an $(a,b,c)$ Type I HBK,
all three exceed $g$, equalling $(g + d)$, $(g + f)$, and $(g + e)$,
for twists between $(A_{H}, B_{H})$, $(B_{H}, C_{H})$ and $(C_{H},
A_{H})$ respectively. For the first such twisting, switch the
L-indices while retaining the U-indices, and consider the product of
$(b, {\bf G} + g + f)$ and $(a, {\bf G} + g + e)$.  Write these
index-pairs in left-to-right order, with the pair starting with $a$
on the bottom, and further mark the right-hand term on the top with
a "+" (since the unit with this index will have \textit{opposite}
sign from that of the other zigzag's corresponding term, hence their
twist product will have the \textit{same} sign in the same
location).  Symbolic multiplication shows combining the four
term-by-term products cannot sum to zero, as follows:

\smallskip
\begin{center}
$+ (b) + ({\bf G} + g + f)$

\underline{$\;\;+ (a)   \;+\;  ({\bf G} + g + e)\;$}
\smallskip

$- ({\bf G} + g + s) {  } + (c)$

\underline{$ + (c) {    } + ({\bf G} + g + s)$}

\smallskip
$2c: \;\;  NOT \; ZERO$\end{center}

\smallskip

Since $(e,f,c)$ is CPO for the L-index trip along the bottom side of
the triangle, and $(a,b,c)$ is likewise the Z-trip in CPO, the upper
right and lower left results are clearly same-signed copies of $c$.
(Since the top right product entails \textit{two} Rule 2
eliminations, there is no resulting sign-change effecting $ef$.) The
other two results, meanwhile, both prepend $({\bf G} + g)$ with no
effect on sign; signs on $s$ are opposite, however, because one
results from multiplying terms in $v \cdot z$ order, while the other
multiplication has form $z \cdot v$.

By symmetry, parallel results (of $2a$ and $2b$) are obtained for
the twists of $B_{H}C_{H}$ and $C_{H}A_{H}$ respectively.  As these
exhaust the nontrivial possibilities for products between L- and U-
index pairings in $O$, our initial claim proves true.

A similar symmetry argument lets us prove our case for all Type I
trefoil HBK's by proving it for any one:  recall that each of these
differs from Type I by having reversed flows along the strut and
side containing the single L-index shared with the zigzag in the
pre-exploded box-kite.  It will suffice, then, to consider the case
of the $(a,d,e)$ HBK derived from the same pre-explosion box-kite.
Its zigzag assessors $A_{H}$, $B_{H}$, $C_{H}$ then read $(a, {\bf
G} + g + f)$; $(d, {\bf G} + g + c)$; $(e, {\bf G} + g + b)$. Twists
among these three all take one to box-kites with strut constants
exceeding $g$, with similar effect.  Since edge-signs are not all
the same (two are positive, one negative), twist-product signs will
also lack uniformity, so we append the binary variable $\sigma$ to
the top right term in the multiplication's setup, and test that
neither possible value for it can lead to a final summation of zero.

Twisting $A_{H}$ and $B_{H}$ (with Fano-plane arrow curving between
their L-indices in the zigzag's circle and beneath $d_{H} = g + b$)
yields the following symbolic arithmetic:

\smallskip
\begin{center}
$+ (d) + \sigma \cdot ({\bf G} + g + a)$

\underline{$\;\;+ (a)   \;+\;  ({\bf G} + g + c)\;$}
\smallskip

$ - ({\bf G} + g + s) \;+\; \sigma \cdot e$

\underline{$ - (e) {    } - \sigma \cdot ({\bf G} + g + s)$}

\smallskip
$2e \; or \; 2({\bf G} + g + s):  NOT \; ZERO$\end{center}

\smallskip

By symmetry once more, we assert a ``not zero'' result for
twist-product zero-divisor testing for $B_{H}C_{H}$ and
$C_{H}A_{H}$. This completes the proof of \textit{universality}: all
zigzags of all Type I HBK's have cowbird's nests, and thereby can
``hold eggs.''

Only the zigzags, however, in such spandrels can do so.  For while
the strut constants of box-kites being twisted to all exceed $g$ in
the above treatments, the zigzag's L-indices all are \textit{less}
than $g$, and two of these will be strut constants for box-kites
being twisted to for the trefoil HBK's.  When the source box-kite
resides in the sedenions, this clearly means these two {\bf S}
values being twisted to designate ET's with no empty spaces save the
long diagonals; hence, all products of ZD's within them that aren't
trivially excluded (self-products and strut-constants) will be zero.
But symbolic calculations like those just considered make it clear
that this low-$N$ situation generalizes completely.

Consider, for instance, the $(a,b,c)$ HBK's own $(a,d,e)$ sail,
$(a_{H}, d_{H}, e_{H}) = a, g+d, g+e$.  By Rule 2, the last two
terms must be placed in reverse order when their trip is written
CPO, so we have $(a, g+e, g+d) = (a_{H}, e_{H}, d_{H})$:  the side
and strut containing $a$ are both reversed in this HBK.  U-indices
are $(G + g + f)$, $(G + c)$, and $(G + b)$. (Note the lack of a `g'
term in both $D_{H}$ and $E_{H}$.) The $a_{H} d_{H}$ twist implies a
box-kite with strut constant $a < g$, leading to this symbolic
arithmetic:

\smallskip

\begin{center}
$+ (a) + \sigma \cdot ({\bf G} + c)$

\underline{$\;\;+ (g + d)   \;+\;  ({\bf G} + g + f)\;$}

$ + ({\bf G} + g + s) \;+\; \sigma \cdot (g + e)$

\underline{$ + (g + e) {    } + \sigma \cdot ({\bf G} + g + s)$}

\smallskip
$( \sigma = -1) \Leftrightarrow 0 \;\; \therefore \; NO \; EGG$
\end{center}

\smallskip
The obvious symmetry argument holds for the other two trefoils'
twist products involving the assessor the reversed side shares with
the zigzag.  Our claim of ``no egg'' must hold, then, for all
trefoil sails within the $(a,b,c)$ HDK.  But as the trefoil-based
HDK's also have (2 of 3) zigzag twist products in box-kites with
strut constants $< g$, mutually zero-dividing ZD's must be contained
in these, hence ``no egg.''  Hence, the Type I HBK eggs are all to
be found nestled in zigzag sails only, proving our claim of
\textit{habitat locality}.

Finally, what of Type II box-kites, which (aside, of course, from
HBK's) we know by Lemma 7 are the only other kind that exist? Recall
from Lemma 6 how they derive from Type I's:  to test their HBK's,
we'll need to consider 3, not 2, high-bit appendings, which we'll
label, by increasing size, $g$, ${\bf G}$, and $\Gamma$.  This makes
the bookkeeping more convoluted, but the proofs are not more
difficult in principle.  This is where the notion of oriented lines
in two different Fano presentations being \textit{flowmorphic},
introduced toward the end of the Cowbird Theorem's statement, must
be made concrete.

As Type I spandrels require only one explosion, symbolic expressions
describing oriented triplets of nodes within them entail only the
use of $g$ in addition to the standard seven letters.  Type II
spandrels, however, require two PL operations, the first to convert
the Type I into a Type II. Hence, symbolic expressions within them
entail ${\bf G}$ as well as $g,$ and so will frequently require
applying Rule 2 twice.  This is best explained through examples. For
Type I, we stick with our running for-instance, exemplified in Table
2:  the explosion of an ${\bf S}=1$ sedenion box-kite with zigzag
L-indices $(3, 6, 5)$. By XOR with {\bf S}, we know that, for the
latter source, the strut opposite $f, e, d$ L-indices are $2, 7, 4$
in that order, with corresponding U-indices found by adding {\bf G}
to the strut opposites:  hence, $A = {\bf G} + f = 10$, $B = {\bf G}
+ e = 15$, $C = {\bf G} + d = 12$, and so on.  Now, explode the (A,
B, C) sail into its own HBK by PL, and we get the leftmost "ABC"
column entry in the top row of the table, presented on the page
facing the statement of the Cowbird Theorem: the old $(a , A) = (3,
10)$ is split into the new strut-opposite $(a, f)$ pairing, and
their U-indices are now the new ${\bf G} + f = 16 + 10 = 26$, with
the rest filled out as shown in the table.  The next three columns
do the same for the trefoil sails' assessors as indicated.  The row
immediately below then spells out the contents of the F, E, D
assessors associated with the A, B, C's above.

Now, let's first re-write the $(A,B,C)$ HBK with $(a,b,c)$
\textit{sail} in the zigzag, both literally and symbolically.
(\textit{Re-}write, since this is what is given by default upon
completing the usual G-loading by corners-and-center PL.)  Then, use
DX($ade; D)$ -- D, not H, since all 3 trefoil L-index sets are
\textit{reversed} in this HBK -- to swap zigzag and $a$-sharing
trefoil while preserving both their orientations, so that we can
compare these two different Fano presentations for this same HBK.

Next, do the same for the Type II built by G-loading (then swapping,
to retain original orientations) $d$ and $e$, and $b$ and $c$, per
the move from the left to right diagrams in Figure 4.  Use $g$ in
the symbolic expressions of this first step, then explode into an
HBK with $2g = {\bf G}$, eventually deriving associated U-indices
via adding of $4g = 2{\bf G} = \Gamma$ to L-index strut-opposites.
Likewise, write the Type II $(A,B,C)$ HBK in two different Fano
presentations, with $(a,b,c)$ and $(a,d,e)$ in their respective
zigzags.

Shorthand the $(a,b,c)$-, then $(d,e,f)$- as-zigzag presentations of
the Type I $(A,B,C)$ HBK as I:ABC $\looparrowright$ $abc$ and I:ABC
$\looparrowright$ $def$ in that order, and write the analogous Type
II's identically, but with II preceding the colon.  If one converts
all symbolic expressions at nodes into pure graphical elements
(oriented arrows), two deep surprises are revealed:

\smallskip
\begin{center}

I:ABC $\looparrowright$ $abc \asymp$ II:ABC $\looparrowright$ $ade$

I:ABC $\looparrowright$ $ade \asymp$ II:ABC $\looparrowright$ $abc$

\end{center}
\smallskip

Here, we use $\asymp$ to mean \textit{is flowmorphic to}:  their
oriented Fano graphs, with nodes labeled only with the minimal
symbolic elements of the set $(a,b,c,d,e,f,s)$, are identical.

Proceed in the same fashion, employing all four columns of literal
indices in the aforesaid table:  taking care to match letters and
positions (so that the $b$ L-index of the third $(f,d,b)$ trefoil is
written in the second position of the zigzag, the $c$ L-index of the
fourth $(f,c,e)$ is put in the third zigzag slot, etc.), the
flowmorphic correspondences continue, in the same exact manner.

Consider the 3 T-bar graphs mentioned earlier, which present as
successive rotations clockwise through $120^{o}$:  two reversed
lines, the T's crossbar and stem, share the $a$ node in I:ADE
$\looparrowright$ $abc$, the $b$ node in I:DBF $\looparrowright$
$abc$, and $c$ in I:EFC $\looparrowright$ $abc$.  These correspond
to Type II's analogous graphs in the $(a,d,e)$ presentations of the
ADE, DBF, and EFC HBK's respectively (with the $(a,d,e)$ in the Type
II ABC being flowmorphic to the Type I HBK's
all-lines-reversed-but-the-zigzag explosion graph).

Ditto, for the $(a,d,e)$ row of presentations in Type I's HBK's and
the $(a,b,c)$ row in Type II's.  For all, the zigzag circular line
is always clockwise (guaranteed by judicious choice of H or D in the
DX).  Viewed left (the ABC HBK) to right (the EFC), we first
encounter the \textit{Pup Tent} graph, so-called because it has all
3 sides reversed, and only the strut \textit{not} reversed in the
pathion Type II box-kite source being reversed here, suggesting the
vertical zipper in a pup-tent's triangular entryway. Then comes the
\textit{Swallow's Tail}, with all non-zigzag lines reversed, except
the two sides with midpoints not $a$.  Third is the two-tined
\textit{Shrimp Fork}, where the only reversals are the line-pair
pointing out and away from $e$ and not including the third ray from
$e$'s angle containing $a$. Finally, in the EFC HBK, the two
reversals of the \textit{Switchblade} comprise the side on the right
leading into, and the strut leading out from, the node $d$ at the
top.

For the third and fourth rows we can construct beneath these same
column heads, DX's bringing properly oriented $(d,b,f)$ and
$(e,f,c)$ L-trips into the zigzag, display T-bar-style rotations of
only these four graphs (Pup-tents always and only in the ABC column,
the others shuffled around in the other columns). This exact
correspondence between Type I and Type II HBK graphs makes it clear
that Type II's second row of graphs, and Type I's first, comprise
the complete set of cowbird's nests \dots thereby completing our
proof.

$\blacksquare$

\medskip

\section{Cowbirding, Bricolage, and Future Directions}

Cowbirds famously lay their eggs in other birds' nests.  As a verb,
``cowbirding'' was how some object-oriented programmers at the old
Lotus Development Corporation described stuffing methods or
structures in abandoned object slots, when creating new ones was
inconvenient or disallowed.  Our cowbird's nests permit infiltration
from outside the current index-system context, directly from another
such: indefinitely many octonion copies, one per spandrel HBK, means
innumerable sites from which to restart CDP.  Map indices of units
in a given nest $(3, 6, 5, 26, 31, 28$, plus ${\bf X} = 25$ and $0$
for Reals in our running example) to the usual ``starter kit'' (of
$0$ through $7$, with ${\bf X}$ mapping most readily to $4$ in the
center, per the L-index Fano of the sedenions' ${\bf S} = 4$
box-kite). Or map them to any other ``kit'' that seemed convenient,
then back again after ``digressing.''  Chomsky's context-sensitive
grammars (as opposed to the context-\textit{free} typical of
programming languages) are clearly implicated -- suggesting
algorithmic opportunities exceeding the built-in givens of our ``new
kind of Number Theory.''

``Cowbirding'' as described here is synonymous with a term made
famous by L\'evi-Strauss in \textit{The Savage Mind}, then
disseminated by his colleague Fran\c{c}ois Jacob [15] among
evolutionary biologists. \textit{Bricolage}, per the anthropologist,
is what a rural jack-of-all-trades or ``Mr. Fixit'' (translated as a
``tinkerer'' in Jacob's piece) performs.  Like ``the significant
images of myth,''

\begin{quote}
the materials of the bricoleur are elements which can be defined by
two criteria:  they have \textit{had a use}, as words in a piece of
discourse which mythical thought `detaches' in the same way as a
bricoleur, in the course of repairing them, detaches the cogwheels
of an old alarm clock; and \textit{they can be used again} either
for the same purpose or for a different one if they are at all
diverted from their previous function. [16, p. 35]
\end{quote}

Like MacGyver, the Swiss army knife and duct-tape-toting protagonist
of the eponymous TV series, the bricoleur, bereft of a specialized
collection of high-tech tools, employs odds and ends he finds at
hand, solving seemingly unrelated problems of the moment in
unconventional ways. In the words of the authors of a highly
influential tract on cognitive science, this is one among many ways
of describing evolution as driven by suboptimal solution-finding, or
``satisficing,'' wherein selection

\begin{quote}
operates as a broad survival filter that admits any structure that
has sufficient integrity to persist.  Given this point of view, the
focus of analysis is no longer on traits but rather on organismic
patterns via their life history.  [This] is evolution as
\textit{bricolage}, the putting together of parts and items in
complicated arrays, not because they fulfill some ideal design but
simply because they are possible.  Here the evolutionary problem is
no longer how to force a precise trajectory by the requirements of
optimal fitness; it is, rather, how to prune the multiplicity of
viable trajectories that exist at any given point. [17, p. 196]
\end{quote}

What's needed:  a grasp of dynamic processes that drive or enable
our flip-books, balloon-rides, explosions, and so forth. Elsewhere
(see pp. 139-40 of [1], and Sections 2 and 3 of [14]), we have
deployed ZD tools to represent key objects of semiotics.  For
example, we used the correspondence between the 4-unit pattern of
strut-opposite assessors and Jean Petitot's 4-control Butterfly
Catastrophe rendering of Algirdas Greimas' ``Semiotic Square.''[15]
But Petitot also provides a more complex, Double Cusp Catastrophe,
model of the primary tool used by L\'evi-Strauss: the ``canonical
law of myths.''[16]  As the third section of this monograph's first
draft sketches out at length[18], sequel studies[19] will deploy a
sort of Catastrophic representation theory based on ZD's. Beginning
with the correspondence between the ``local'' level of our
strut-opposite quartets of assessor indices and the Semiotic Square,
we proceed into more rarefied air, where explosion necessitates
connecting with the ``global'' network-trafficking the Canonical Law
would regulate. (Leading question: Will the future architecture of
the Web-replacing Grid recapitulate that of L\'evi-Strauss's ``web
that knows no weaver'' made of myths?)

To briefly review our earlier work relating ZD patterns to Greimas'
Square, place $A$ and $a$, say, on a square's top corners, and strut
opposite units $F$ and $f$ on the corners diagonally opposite.  The
two ends of each line are connected by XOR-ing with {\bf X}, the
diagonals by {\bf G}, and the verticals by {\bf S}. Instantiations
of this ``atom of meaningfulness'' abound in Greimas and Petitot's
works, but one telling example implicitly demonstrates its
difference from Boolean binary logic.

Across the horizontals at the top and bottom of the box, write
``True'' and ``False''; along left and right verticals, put
``Secret'' and ``Lie.''  Label the nodes at top left and bottom
right ``Being'' and ``Non-Being,'' and refer to the diagonal as the
schema for \textit{immanence}; those at top and bottom of the other
diagonal, regulating \textit{manifestation}, are labeled ``Seeming''
and ``Non-Seeming'' respectively.  (For a detailed discussion, see
Section 3.7 of [15]). This provides a framework for contemplating
verediction, which plays a key role in the contractual component of
narratives.  The exposure of the villain transforms ``Lie'' into
``False'' at the turning point in countless fairy tales, where the
threat to the true order of things is finally rejected (typically,
at the last possible moment).  In stories like \textit{Cinderella},
the narrative is propelled by the inevitability of transforming the
``Secret'' into the ``True'':  it is the possessor of the secret of
the glass slipper, not one of her evil step-sisters, who rightly
wins Prince Charming's heart.

The three kinds of lines relate to Roman Jakobson's three kinds of
``binary opposition'' in his groundbreaking studies of phonemics, at
the basis of all later structuralist set-ups, including that of
L\'evi-Strauss.  The diagonals indicate Jakobson's ``privation'':
e.g., the plosive 'p' differs from 'b' solely by its absence of
voicing -- and indicate, for \textit{us}, where the singleton
high-bit indicated by {\bf G} is XOR'd with the index of a lower- or
upper- case letter, thereby connecting L- to U- units in
strut-opposite assessors.

The horizontals are sites of \textit{contrariety} -- a 2-control
competition between 2 warring parties in Catastrophe terms (or a
pair of ``sememes'' forced into relationship from the semiotician's
vantage).  They generate the synchronous (\verb|/|) and
anti-synchronous (\verb|\|) diagonals in the assessor planes in our
model.  Verticals are linked with \textit{implication}.  Per last
paragraph's examples, transforming competitive dynamics on
horizontals, into orders of implication along verticals, opens the
door to higher-order models: conversion of horizontals into
verticals in this sense is exactly what our \textit{explosion
process} effects.

Here we can underwrite the full workings of ``spandrel thinking,''
as Gould and Lewontin explain it. For spandrels exist not only on
postage stamps, but in the quartet of curved triangles formed where
dome-supporting arches cross in front of cathedral naves. Spandrels
became favorite sites for mosaic and painterly expression -- so much
so, that one who was architecturally naive might think the archways'
intersection pattern was concocted to facilitate their production.
But in fact, they are the happy side-effect of the architecture; the
evolution of architectural design selected for crossed arches, not
the spandrels that rode on their coattails. Gould and Lewontin's
point:  many evolutionary arguments assume selection pressures are
at work evolving spandrel-like attributes -- or, in Greimas' argot,
that presuppositions (the Square's verticals) must sometimes fight
for survival (along horizontals).

Such cart-before-horse flipflops are endemic in any explanatory
enterprise.  What we claim here is our toolkit suffices to model
conundra of this sort... and allow for contexts wherein spandrels,
by cowbird logic, become sites for future adaptations (hence,
selection pressures) in their own right.[20]

\section*{References}

\begin{description}

\item \verb|[1]| Robert P. C. de Marrais, ``Placeholder
Substructures I:  The Road from NKS to Scale-Free Networks is Paved
with Zero-Divisors,'' \textit{Complex Systems, 17 (2007)},125-142;
arXiv:math.RA/0703745.

\item \verb|[2]| Robert P. C. de Marrais, ``Placeholder
Substructures II:  Meta-Fractals, Made of Box-Kites, Fill
Infinite-Dimensional Skies,'' arXiv:0704.0026 [math.RA]

\item \verb|[3]| Robert P. C. de Marrais, ``Placeholder
Substructures III:  A Bit-String-Driven `Recipe Theory'...,''
arXiv:0704.0112 [math.RA]

\item \verb|[4]| Fran\c{c}ois Jacob, ''Evolution and Tinkering,''
\textit{Science}, 196 (1977), 1161-1166.

\item \verb|[5]| Claude L\'evi-Strauss, \textit{The Savage Mind}
(University of Chicago Press, Chicago and London, 1969; French
original, 1962)

\item \verb|[6]| Robert P. C. de Marrais, ``The 42 Assessors and
the Box-Kites They Fly,'' arXiv:math.GM/0011260

\item \verb|[7]| R. Guillermo Moreno, ``The Zero Divisors of the
Cayley-Dickson Algebras over the Real Numbers,'' \textit{Boletin
Sociedad Matematica Mexicana (3)}, 4, 1 (1998), 13-28; preprint
available online at arXiv:q-alg/9710013

\item \verb|[8]| Sir Tim Berners-Lee, ``The Two Magics of Web
Science,'' WWW Banff keynote address, May 9, 2007; video available
online at www.w3.org/2007/ \linebreak Talks/0509-www-keynote-tbl

\item \verb|[9]| Robert P. C. de Marrais, ``Placeholder
Substructures," wolframscience.com/ \linebreak
conference/2006/presentations/materials/demarrais.ppt (Powerpoint
show)

\item \verb|[10]| I. L. Kantor and A. S. Solodovnikov, \textit{Hypercomplex
Numbers:  An Elementary Introduction to Algebras} (Springer-Verlag,
New York, 1989)

\item \verb|[11]| Stephen Jay Gould and Richard C. Lewontin, ``The
Spandrels of San Marco and the Panglossian Paradigm:  A Critique of
the Adaptationist Programme,'' \textit{Proceedings of the Royal
Society of London}, Series B, Vol. 205, No. 1161 (1979), 581-598.

\item \verb|[12]| Robert P. C. de Marrais, ``Flying Higher Than A
Box-Kite,'' arXiv:math.RA/ \newline 0207003

\item \verb|[13]| Benoit Mandelbrot, \textit{The Fractal Geometry of Nature} (W.
H. Freeman and Company, San Francisco, 1983)

\item \verb|[14]| Robert P. C. de Marrais, ``Presto! Digitization,'' arXiv:math.RA/0603281

\item \verb|[15]| Jean Petitot, \textit{Morphogenesis of Meaning},
(Peter Lang, New York and Bern, 2004; French original, 1985).

\item \verb|[16]| Jean Petitot, ``A Morphodynamical Schematization of
the Canonical Formula for Myths,'' in \textit{ The Double Twist:
From Ethnography to Morphodynamics }, edited by Pierre Maranda
(University of Toronto Press, Toronto and Buffalo, 2001, pp.
267-311).

\item \verb|[17]| Francisco J. Varela, Evan Thompson, and Eleanor
Rosch, \textit{The Embodied Mind:  Cognitive Science and Human
Experience} (The MIT Press, Cambridge MA and London, 1991)

\item \verb|[18]| Robert P. C. de Marrais, ``Voyage By Catamaran: Long-Distance Network Navigation, from Myth Logic
  to the Semantic Web, Can Be Effected by Infinite-Dimensional Zero-Divisor
  Ensembles'', arXiv:math.GM/0804.3416 V1 [math.GM].

\item \verb|[19]| Robert P. C. de Marrais, "Natural Numbers, Natural Language:  Architecting the Semantic
Web," invited presentation at 4th International Conference on
Natural Computing (ICNC '08) at Shandong University, Jilin, China,
October 18-20, 2008.

\item \verb|[20]| Adam S. Wilkins, ``Between `design' and
`bricolage':  Genetic networks, levels of selection, and adaptive
evolution,'' PNAS, 104, suppl. 1, May 15, 2007;
http://pnas.org/cgi/doi/ 10/1073/pnas.0701044104

\end{description}

\end{document}